\documentclass{article}
\usepackage[utf8]{inputenc}

\title{Trivial coloring of Cartesian product of graphs}
\author{Tamás Csernák}
\date{January 2023}

\usepackage{amsthm}
\usepackage{amsfonts}

\usepackage{mathrsfs}
\usepackage{enumerate}
\usepackage{amsmath}
\usepackage{amssymb}

\usepackage{todonotes}

\makeatletter
\@namedef{subjclassname@2020}{%
  \textup{2020} Mathematics Subject Classification}
\makeatother
\newtheorem{theorem}{Theorem}[section]        

\newtheorem{lemma}[theorem]{Lemma}
\newtheorem{observation}[theorem]{Observation}
\newtheorem{problem}[theorem]{Problem}

\theoremstyle{definition}
\newtheorem{definition}[theorem]{Definition}

\begin{document}

\maketitle

\begin{abstract}
  A coloring of a direct product of graphs is said to be {\em trivial} iff it is induced by some coloring  of a factor of the product. A graph $G$ is  {\em trivially power colorable} iff 
  every coloring of a finite power of $G$ with $\chi(G)$-many colors is trivial. 
  Greenwell and Lovász proved that the finite complete graphs $K_n$ for $n\ge 3$ 
  are trivially power colorable. Generalizing their result we define a much wider  class of trivially power-colorable graphs: 
  if $G$ is a finite, connected  graph with $\chi(G)\ge 3$ and 
  every vertex of $G$ is in a clique of size $\chi(G)$, then $G$ is  trivially power-colorable.
   As an application of this result,   we give a complete characterization of 
   trivially power-colorable  cographs.

  Finally, we give a structural description of the  colorings of  infinite powers of trivially power-colorable finite graphs.
\end{abstract}


\section{Introduction}



In this paper, we define  a graph as a pair $G=(V(G),E(G))$, where $V(G)$ denotes the  set of vertices and $E(G)\subseteq [V(G)]^{2}$ represents  the set of edges. All graphs we discussed in this work   are undirected and simple.

For any (finite or infinite) index set $I$, and graphs $(G_{i})_{i\in I}$, we define their product as $G=\times_{i\in I}{G_{i}}$, where $V(G)=\times_{i\in I}{V(G_{i})}$, and for $u,v\in V(G)$ with $u=(u_{i})_{i\in I}$ and $v=(v_{i})_{i\in I}$, we have that $u,v\in E(G)$ if and only if $u_{i}v_{i}\in E(G_{i})$ for all $i\in I$. The product of two graphs, $G$ and $H$ can be simply denoted by $G\times H$. If for all $i\in I$, we have $G_{i}=G$, then we can use the notion $G^{I}=\times_{i\in I}{G}$. 

Furthermore,  for any graph, we denote its chromatic number by $\chi (G)$.  
For any $k$ and a finite graph $G$, let $P(G,k)$ be 
the  number of $k$-colorings of $G$.


Finally,  let $c(G)$ denote  the number of connected components of a graph $G$.

\section{Trivially colorable graphs}







For any two graphs $G$ and $H$, it is clear that $\chi (G\times H)\le min(\chi (G),\chi(H))$.
For long time it was conjectured, that for all pairs of finite graphs, we have $\chi (G\times H)= min(\chi (G),\chi(H))$. However, as shown in \cite{YS}, this was recently proven false. 
Nevertheless,  the special case $\chi (G^{n})=\chi (G)$ clearly holds for all graphs $G$. 


\begin{definition} \label{trivial_coloring}
 Let $I$ be an index set and for each $i\in I$, let $G_{i}$ be a graph.   Write $G=\times_{i\in I}{G_{i}}$. The coloring $\Phi$ of $G$   is {\em trivial} if there exists some $i^*\in I$ and a coloring $\phi $ of $G_{i^{*}}$, such that for each  $v=(v_i)_{i\in I}\in V(G)$, we have $\Phi (v)=\phi (v_{i^*})$. The product graph $G$ is said to be {\em trivially  colorable}, if all the $\chi(G)$-colorings of $G$ are trivial. 
 
 A graph $H$ is {\em trivially power-colorable}, if for all $n\in \omega $, the product $H^{n}$ is trivially colorable. 
\end{definition}



The solution of the $n$-switch problem,  (see Greenwell and Lovász \cite{GLL})  can be stated as follows: {\em the complete graph $K_{k}$ is trivially power-colorable for all $k\ge 3$.}

\begin{definition}\label{tight}

  Let $G$ be a graph, and let $\phi$ be a $k$-coloring of $G$. We say that $\phi$ is tight if, for every vertex $v \in V(G)$ and every color $c \neq \phi(v)$, there is a neighbor $u$ of $v$ such that $\phi(u) = c$. In other words, $\phi$ is tight if no vertex's color can be changed without violating the proper coloring condition.

  $G$ is {\em tight} if all of its $\chi(G)$-colorings are tight.
\end{definition}

Greenwell and Lovász has also shown, that if $H$ is tight with $\chi (H)=k$ for some $k\ge 3$, then  $K_{k}\times H$ is trivially colorable (see \cite[Theorem 1]{GLL}).  

The primary aim  of this paper is to characterize  a broad class of the finite graphs that are trivially power-colorable, or at least to provide  some necessary and some sufficient conditions for this property. 

In the final section,  we examine  the colorings of  graphs s$G^{\lambda }$, where $G$ is a finite trivially power-colorable graph and ${\lambda}$ is an infinite cardinal. 

\subsubsection*{Preliminary observations}

First,  we show that trivially power-colorable graphs must be connected. Some of these statements are already known, but for the completeness, I will present short proofs.

\begin{lemma} \label{product_component}
 If $G$ is a finite graph, then for all $n$, we have $c(G^{n})\ge c(G)^{n}$.
\end{lemma}


\begin{proof}
 Let us denote for any graph $H$ and $u,v\in V(H)$, the equivalence relation $u\sim v$, that $u,v$ are in the same component in $H$. If $u,v\in V(G^{n})$ and $uv\in E(G^{n})$, then $u_{i}v_{i}\in E(G)$ for all $1\le i\le n$, so clearly $u_{i}\sim v_{i}$. Since $\sim $ is an equivalence relation through the edges, we have, that if for $u,v\in V(G^{n})$ with $u\sim v$, then $u_{i} \sim v_{i}$ for all $i$. Choose $w_{1},...,w_{c(G)}\in V(G)$, such that they are all in different components, and let $U=\{ (w_{j_{1}},...,w_{j_{n}}): 1\le j_{1},...,j_{n}\le c(G)\} \subseteq V(G^{n}) $. Clearly $|U|=c(G)^{n}$ and if $u,v\in U,u\neq v$, then there is some $i$, such that $u_{i}\neq v_{i}$, and since, they are chosen from different components, we have $u_{i}\not\sim v_{i}$, so $u\not\sim v$. All elements of $U$ are in pairwise different component of $G^{n}$, so $c(G^{n})\ge |U|=c(G)^{n}$. 
\end{proof}

\begin{lemma}\label{chromatic_polynomial_component}
 If $G$ is a finite graph and $\chi (G)\le k$, then $P(G,k)\ge k^{c(G)}$.
\end{lemma}

\begin{proof}
 Let $r=c(G)$, and $S_{1},...,S_{r}\subseteq V(G)$ be the connected components. Then all induced subgraphs on $S_{i}$ are $k$-colorable. By the permutation of colors used, we can easily see that $P(G[S_{j}],k)\ge k$ for all $j$. 
 All the combinations of $k$ colorings of components gives a $k$ coloring of $G$, so 
 $P(G,k)=\prod^r_{j=1}P(G[S_{j}],k)\ge k^{r}=k^{c(G)}.$
\end{proof}

\begin{theorem}\label{trivial_connected}
 If $G$ is a trivially power-colorable finite graph, then $G$ is connected.
\end{theorem}

\begin{proof}
Write $\chi (G)=k$.  Since  $G$ is trivially power-colorable, for all $k$-colorings $\Phi $ of $G^{n}$ there is $1\le i\le n$ and a $k$-coloring, $\phi $ of $G$, such that $\Phi (v)=\phi (v_{i}) $ for all $v\in V(G^{n})$. This $i$ can be chosen $n$ different ways and $\phi $ can be chosen $P(G,k)$ different ways, so we have $nP(G,k)\ge P(G^{n},k)$. 

On the other hand, using Lemma \ref{chromatic_polynomial_component} for $G^{n}$, and Lemma 
\ref{product_component} for $G$, we get $P(G^{n},k)\ge k^{c(G^{n})}\ge k^{c(G)^{n}}$.
Putting together, we obtain $nP(G,k)\ge  k^{c(G)^{n}}$.
If $G$ were disconnected, then $c(G)\ge 2$, the right side would be a doubly exponential function of $n$, while the left side would be linear, so for sufficiently large $n$, it would cause a contradiction. Thus $G$ must be connected.
\end{proof}

Now, we have seen, that $G$ must be connected. However it is also possible, that even if $G$ is connected, some of their powers are not.

\begin{lemma}\label{bipartite_disconnected}
If $G$ is a bipartite graph, then $G^{2}$ is not connected.\cite[Theorem 1]{PMW}
\end{lemma}

\begin{proof}
 Let $\phi $ be a 2-coloring of $G$ and $U_{1}=\{ (v_{1},v_{2}):\phi (v_{1})=\phi (v_{2})\} \subseteq V(G^{2})$, $U_{2}=\{ (v_{1},v_{2}):\phi (v_{1})\neq \phi (v_{2})\} \subseteq V(G^{2})$. Then if $uv\in E(G^{2})$, then $\phi (u_{1})\neq \phi (v_{1})$, $\phi (u_{2})\neq \phi (v_{2})$, so if $\phi (u_{1})= \phi (u_{2})$, then $\phi (v_{1})= \phi (v_{2})$, and if $\phi (u_{1})\neq \phi (u_{2})$, then $\phi (v_{1})\neq \phi (v_{2})$. So both are in $U_{1}$ or both are in $U_{2}$, hence these sets form a partition with no edges, thus $G^{2}$ is disconnected.
\end{proof}

\begin{theorem}\label{trivial_chromatic_3}
If $G$ is a trivially power-colorable non empty finite graph, then $\chi (G)\ge 3$.
\end{theorem}

\begin{proof}
 Suppose for contradiction, that $G$ is bipartite, then by Lemma \ref{bipartite_disconnected}, $G^{2}$ is disconnected. Since $G$ is trivially power-colorable, for all $n$, the colorings of $(G^{2})^{n}=G^{2n}$ are trivial in respect of $G$, so they are clearly trivial in respect of $G^{2}$, thus $G^{2}$ is trivially colorable, that contradicts Theorem \ref{trivial_connected}.
\end{proof}

We have seen, that if $G$ is trivially power-colorable, then $G$ must be connected and $\chi (G)\ge 3$. These are necessary conditions.

On the other hand if $G$ is connected and   $\chi (G)\ge 3$, then by \cite[Theorem 1]{PMW} all  
the powers $G^{n}$ are connected. Next, we will see an other necessary condition for $G$ to be trivially power-colorable.


\begin{theorem}\label{trivial_tight}
 If $G^{2}$ is a trivially power-colorable for a finite graph $G$, then $G$ is tight. Hence, trivially power-colorable graphs must be tight.
\end{theorem}


\begin{proof} Write $k=\chi(G)$.
 Suppose for the sake of contradiction, that $G$ is not tight, and let $\phi,\phi'$ be two $k$-colorings of $G$ differ at one point, say $w\in V(G)$. 
 
 Now, we construct a non-trivial $k$-coloring $\Phi $ of $G^{2}$ as follows:
 \begin{displaymath}
 {\Phi(u,v)}=\left\{\begin{array}{ll}
 {\phi'(v)}&\text{if $(u,v)=(w,w)$, }\\
 {\phi(v)}&\text{otherwise.
 }\\
 \end{array}\right.
 \end{displaymath}

 
 First, we need to verify that  $\Phi$ is a good coloring. Let $(u,v)(u',v')\in E(G^{2})$, so 
 $uu',vv''\in E(G)$. 
 If $(u,v)\ne (w,w)$ and $(u',v')\ne (w,w)$, then 
 \begin{displaymath}
 \Phi(u,v)=\phi(v)\ne \phi(v')=\Phi(u',v').
 \end{displaymath}
 If $(u,v)=(w,w)$, then $v'\ne w$ and so 
 \begin{displaymath}
  \Phi(u,v)=\Phi(w,w)=\phi'(w)\ne \phi'(v')=\phi(v')=\Phi(u',v').
  \end{displaymath}
 The case $(u',v')=(w,w)$ is analogue.


 Finally, we show that $\Phi $ is non-trivial, meaning the color of a vertex is not solely determined by either the first or second coordinate.

Let $vv'\in E(G)$. We can assume that $v\ne w$.
Then $\Psi(v,v)=\phi(v)\ne \phi(v')=\Phi(v,v')$, so 
the first  coordinate does not determine 
 the color of a vertex in the product.

 Since $\Phi(w,w)=\phi'(w)\ne \phi(w)=\Psi(v,w)$, the second coordinate does not determine 
 the color of a vertex in the product. 
 
 
 Hence, $\Phi $ is non-trivial, contradicting the assumption that $G^2$ is trivially colorable. 
\end{proof}

\section{Cliqued graphs}

In this section we give some sufficient conditions on graphs to be  trivially power-colorable.

\begin{definition}\label{cliqued}
 Let $G$ be a graph with chromatic number $\chi (G)=k$. We say that  $G$ is {\em weakly-cliqued} if for every vertex  $v\in V(G)$ there is a clique $X\subset V(G)$  of size $k$ that contains   $v$. The graph $G$ is {\em strongly-cliqued}, if it has no isolated vertices and for every edge   $uv\in E(G)$, there is a clique $X\subset V(G)$ of size $k$ that contains 
 both $u$ and $v$. 
\end{definition}


Clearly, all strongly cliqued graphs are also weakly cliqued. Indeed, if $G$ is strongly cliqued and $v \in V(G)$, we can find some $u \in V(G)$ with $uv \in E(G)$ because $G$ has no isolated vertices, and so we have  clique of size $k$ containing both $u$ and $v$.


\begin{observation}\label{obs:wc2tigth}
A weakly cliqued graph is tight. 
\end{observation}
Indeed, if $G$ is weakly-cliqued graph with $\chi(G)=k$, and $v\in V$, then 
$G$ contains a cliques  $X$ of size $k$ which contains $v$. Thus, if 
$\phi:V(G)\to k$ is a proper coloring of $G$, then 
$\phi(v)$ must be  the only element of $k\setminus \{\phi(u):u\in X\setminus \{v\}\}$. In other words,   $\phi(v)$ is uniquely determined by $\phi \restriction V\setminus \{v\}$.


\medskip

In this section, we will see, that any weakly-cliqued connected graph $G$ with $\chi (G)\ge 3$ is trivially colorable (see Theorem   \ref{weakly_cliqued_power_trivial}). 

\begin{definition}\label{clique_connected}
 If $G$ is a weakly-cliqued graph, with $\chi (G)=k$ and $X,X'\subseteq V(G)$ are cliques of size $k$, then a {\em clique-path} between $X$ and $X'$ is a sequence $X=X_{0},X_{1},...,X_{r}=X'$ of size $k$ cliques, such that $X_{j}\cap X_{j+1}\neq \emptyset $ for $0\le j\le r-1$. $G$ is {\em clique-connected} if for all pair $X,X'$ of size $k$ cliques, there is a clique-path between them. This way we can also define {\em clique-connected components} of $G$. For 2 vertices $u,v\in V(G)$ a clique-path between $u$ and $v$ is a clique-path between some $X,X'$ of size $k$ cliques, with $u\in X,v\in X'$. On the other hand if there is a clique-path between $u,v$, then for any $X,X'$ of size $k$ cliques, with $u\in X,v\in X'$ there is a clique path between $X,X'$. As all vertices are contained in a size $k$ clique, we can define clique-connectedness and clique-connected components on $V(G)$.
\end{definition}

\begin{lemma}\label{strong_cliqued_connected}
 If $G$ is a strongly-cliqued, connected graph, then it is clique-connected.
\end{lemma}

\begin{proof}
Straightforward.
\end{proof}

 First we try to generalize the theorem of Greenwell and Lovász (Theorem 1 in \cite{GLL}).

\begin{lemma}\label{strongly_cliqued_tight_trivial}
 If $G$ is a clique-connected graph, $H$ is a tight graph with $\chi (G)=\chi (H)=k\ge 3$ , then all $k$-colorigs $\Phi $ of $G\times H$ are trivial. 
\end{lemma}

\begin{proof}
 Let $\Phi:G\times H\rightarrow k$ be a $k$-coloring. For any size $k$ clique $X\subseteq V(G)$, let us define $X$ is of 1st type if there is a coloring $\phi _{X}:X\rightarrow k$  , such that for all $u\in X,v\in V(H)$, we have $\Phi(u,v)=\phi_{X}(u)$, and $X$ is of 2nd type if there is a coloring $\psi _{X}:V(H)\rightarrow k$, such that for all $u\in X,v\in V(H)$, we have $\Phi(u,v)=\psi_{X}(v)$. By Theorem 1 in \cite{GLL}, we have that all size $k$ cliques are either 1st type or 2nd type. 
 
 First we will see, that we cannot have 1st and 2nd type cliques at the same time. Suppose for contradiction, that $X\subseteq V(G)$ is of 1st type, and $X'\subseteq V(G)$ is of 2nd type. Since $G$ is clique-connected, there is some clique path $X=X_{0},...,X_{r}=X'$. All of $X_{j}$-s are either 1st type or 2nd type, so there is some $0\le j\le r-1$, such that $X_{j}$ is of 1st type and $X_{j+1}$ is of 2nd type. Since $X_{j}\cap X_{j+1}\neq \emptyset $, we can choose $u\in X_{j}\cap X_{j+1}$. Pick any $v,v'\in V(H)$ with $vv'\in E(H)$. Then we have $\psi _{X_{j+1}}(v)=\Phi (u,v)=\phi _{X_{j}}(u)=\Phi (u,v')=\psi _{X_{j+1}}(v')$ in contradiction with $\psi_{X_{j+1}}$ is a coloring of $H$. Thus either all cliques are of 1st type or all cliques are of 2nd type.
 
 First suppose that all cliques are of 1st type. for any cliques $X,X'\subset V(G)$ with $X\cap X'\neq \emptyset $ and $u\in X\cap X'$, we have $\phi _{X}(u)=\Phi (u,v)=\phi _{X'}(u)$, where $v\in V(H)$ can be any vertex, so the functions $\phi _{X}$ are compatible. Let $\phi=\bigcup_{X\subseteq V(G),clique}{\phi _{X}}$. Then $\phi :V(G)\rightarrow k$ is a function. It is also clear that for any $u\in V(G)$, there is some size $k$ clique $X\subseteq V(G)$  with $u\in X$, so for all $v\in V(H)$, we have $\Phi (u,v)=\phi _{X}(u)=\phi (u)$. We need to show that this is a good coloring. Let $u,u'\in V(G)$ with $uu'\in E(G)$ and pick any $v,v'\in V(H)$ with $vv'\in E(H)$. Then since $(u,v)(u',v')\in E(G\times H)$, we have $\phi (u)=\Phi(u,v)\neq \Phi (u',v')=\phi (u')$, so $\phi $ is a good coloring.
 
 Now suppose that all cliques are of 2nd type. We need to show, that for all $X,X'\subseteq V(G)$ cliques, we have $\psi _{X}=\psi _{X'}$. If $X\cap X'\neq \emptyset $, then let $u\in X\cap X'$. Then for all $v\in V(H)$, we have $\psi _{X}(v)=\Phi (u,v)=\psi _{X'}(v)$, thus $\psi _{X}=\psi _{X'}$. Now let $X,X'\subseteq V(G)$ be arbitrary and let $X=X_{0},...,X_{r}=X'$ be a clique path. Then $\psi _{X}=\psi _{X_{0}}=\psi _{X_{1}}=...=\psi _{X_{r}}=\psi _{X'}$. Then let $\psi =\psi _{X}$ for all cliques. Then for any $u\in V(G)$ and $v\in V(H)$ choose a size $k$ clique $X$ with $u\in X$, and we have $\Phi (u,v)=\psi _{X}(v)=\psi (v)$.

\end{proof}
 
 Next, we generalize this lemma, by changing the clique-connected property to weakly-cliqued.
 

 \begin{theorem}\label{weakly_cliqued_tight_trivial}
 If $G$ is a weakly-cliqued, connected graph, $H$ is a tight graph with $\chi (G)=\chi (H)=k\ge 3$ , then all $k$-colorings $\Phi $ of $G\times H$ are trivial. 
\end{theorem}

\begin{proof}
Let $\Phi:G\times H\rightarrow k$ be a $k$-coloring. 
For any clique-connected component $U\subseteq V(G)$, let us define $U$ is of 1st type if there is a coloring $\phi _{U}:U\rightarrow k$ , such that for all $u\in U,v\in V(H)$, we have $\Phi(u,v)=\phi_{U}(u)$, and $U$ is of 2nd type if there is a coloring $\psi _{U}:V(H)\rightarrow k$, such that for all $u\in U,v\in V(H)$, we have $\Phi(u,v)=\psi_{U}(v)$. By Lemma \ref{strongly_cliqued_tight_trivial}, we have that all clique-connected components are either 1st type or 2nd type. 

We can also define for an $u\in V(G)$, that $u$ is of 1st type or 2nd type, by the unique clique-connected components $U\subseteq V(G)$ with $u\in U$ is of 1st type or 2nd type. We will show that there is no edge between 1st type and 2nd type vertices. Suppose for contradiction, that $u,u'\in V(G)$ $u$ is of 1st type, $u'$ is of 2nd type and $uu'\in E(G)$. Let $U$ be the unique clique-connected component with $u\in U$ and $U'$ be the unique clique-connected component with $u' \in U'$. Then $U$ is of 1st type and $U'$ is of 2nd type. Since $\chi (H)=k$ all of its $k$ colorings must take all values, so there is some $v'\in V(H)$, such that $\psi _{U'}(v')=\phi _{U}(u)$. Pick any $v\in V(H)$ with $vv'\in E(H)$. Then $(u,v)(u',v')\in E(G\times H)$. On the other hand, we have $\Phi(u,v)=\phi _{U}(u)=\psi _{U'}(v')=\Phi (u',v')$, that is a contradiction. Now since $G$ is connected and there is no edge between 1st type and 2nd type vertices, we must have either all vertices are of 1st type or all vertices are of 2nd type, thus all clique-connected components are of 1st type or all clique-connected components are of 2nd type. 

First suppose the case all clique-connected components are of 1st type. Then all the $\phi _{U}$ functions are defined on pairwise disjoint sets, so their union $\phi =\bigcup_{U\subseteq V(G), component}{\phi _{U}}$ is a is a function $\phi :V(G)\rightarrow k$. Clearly if $u\in V(G),v\in V(H)$, we can pick the unique clique-connected component $U\subseteq V(G)$ with $u\in U$, and then we have $\Phi (u,v)=\phi _{U}(u)=\phi (u)$. We need to show that $\phi $ is a good coloring. Let $u,u'\in V(G)$ with $uu'\in E(G)$. Pick any $v,v'\in V(H)$ with $vv'\in E(H)$. Then $(u,v)(u',v')\in E(G\times H)$, so $\phi (u)=\Phi (u,v)\neq \Phi (u',v')=\phi (u')$, thus $\phi $ is a good coloring. 

Now suppose that all clique-connected components are of 2nd type. We need to show, that for any $U,U'\subseteq V(G)$ clique-connected components $\psi _{U}=\psi _{U'}$ holds. First suppose that there is some $u\in U,u'\in U'$ with $uu'\in E(G)$ and suppose for contradiction, that $\psi _{U}\neq \psi _{U'}$. Then there is some $v\in V(H)$ with $\psi _{U}(v)\neq \psi _{U'}(v)$. Since $H$ is tight, we cannot change the color of the vertex $v$ in $\psi _{U'}$ to $\psi _{U}(v)$, so there must be some $v'\in V(H)$ with $vv'\in E(H)$ and $\psi _{U'}(v')=\psi _{U}(v)$. Then $(u,v),(u',v')\in E(G\times H)$. On the other hand, we have $\Phi (u,v)=\psi _{U}(v)=\psi _{U'}(v')=\Phi (u',v')$, that is a contradiction. Now let $U,U'\subseteq V(G)$ be arbitrary clique-connected components and pick $u\in U,u'\in U'$. Since $G$ is connected, there is a path $P$ from $u$ to $u'$. Let $U=U_{0},U_{1},...,U_{r}=U'$ be the list of clique-connected components in the order $P$ passing through,  so for all $0\le j\le r-1$, we have that there is an edge between $U_{j},U_{j+1}$, So $\psi _{U_{j}}=\psi _{U_{j+1}}$. Then $\psi _{U}=\psi _{U_{0}}=\psi _{U_{1}}=...=\psi_{U_{r}}=\psi _{U'}$. Let $\psi =\psi _{U}$ for all clique-connected components. Then for all $u\in V(G),v\in V(H)$, let $U\subseteq V(G)$ be the unique clique-connected component with $u\in U$, we have $\Phi (u,v)=\psi _{U}(v)=\psi (v)$.
\end{proof} 

\begin{theorem}\label{weakly_cliqued_power_trivial}
 If $G$ is a weakly-cliqued, connected graph with $\chi (G)\ge 3$, then $G$ is trivially power-colorable. 
\end{theorem}

\begin{proof}
 Let $k=\chi (G)$. We clearly have that $\chi (G^{n})=k$ for all $n$. By \cite{PMW}, we have that $G^{n}$ is also connected for all $n$. Now we will show, that for all $n$, $G^{n}$ is weakly-cliqued. Let $u\in V(G^{n})$. Then $u=(v_{1},...,v_{n})$, where $v_{j}\in V(G)$ for all $1\le j\le n$. For all $j$ let $X_{j}\subseteq V(G)$ be a clique of size $k$ with $v_{j}\in X_{j}$, and let us list the elements $X_{j}=\{ w_{j,1},...,w_{j,k}\}$ where $w_{j,1}=v_{j}$. For all $1\le i\le k$ let $u_{i}=(w_{1,i},...,w_{n,i})$. Then clearly by definition, we have $u_{1}=u$ and for all $1\le i,i'\le k$ with $i\neq i'$, we have that for all $1\le j\le n$, since $w_{j,i},w_{j,i'}\in X_{j}$ and $w_{j,i}\neq w_{j,i'}$, we have $w_{j,i}w_{j,i'}\in E(G)$, thus $u_{i}u_{i'}\in E(G^{n})$. Then $X=\{ u_{1},...,u_{n}\}$ is a clique of size $k$ in $G^{n}$ containing $u$, so $G^{n}$ is weakly-cliqued.
 
 Now we will show that for all $n$ all colorings of $G^{n}$ are trivial. We prove it by induction on $n$. For $n=1$, it is clear. Suppose it is true for some $n$, and we will prove it for $n+1$. Let $\Phi :V(G^{n+1})\rightarrow k$ be a $k$-coloring of $G^{n+1}=G^{n}\times G$. Since all conditions of Theorem \ref{weakly_cliqued_tight_trivial} hold, by Theorem  \ref{weakly_cliqued_tight_trivial}, we have that either there is some $\tau :V(G^{n})\rightarrow k$ good coloring, such that for all $v\in V(G^{n+1})$, we have $\Phi (v)=\tau(v|_{n})$ or there is a $\psi :V(G)\rightarrow k$, such that for all $v\in V(G^{n+1})$, we have $\Phi (v)=\psi(v_{n+1})$. In the former case by induction there is some $1\le j\le n$, and $\phi :V(G)\rightarrow k$ coloring, such that $\tau(w)=\phi (w_{j})$ for all $w\in V(G^{n})$. But then for any $v\in V(G^{n+1})$ we have $\Phi (v)=\tau (v|_{n})=\phi ((v|_{n})_{j})=\phi (v_{j})$, so $\Phi $ is trivial by coordinate $j$. In the latter case $\Phi $ is clearly trivial by coordinate $n+1$, so the induction step works.
\end{proof}


\section{An application: trivially colorable cographs}

In Section 2, we gave some necessary conditions for a graph to be trivially colorable, and in Section 3, we gave some sufficient conditions. However, these conditions do not characterize trivially colorable graphs in general, as there are many graphs, that meet the necessary conditions but not the sufficient ones. This is not the case, in a special class of graphs, called cographs.

\begin{definition}\label{cograph}
 The {\em cographs} are the finite graphs that can be constructed by the following 3 steps:
\begin{enumerate}[1)]
\item   Any graph with one vertex is a cograph.
 
\item  If $G$ is a cograph, then its complement $\bar{G}$ is a cograph.

\item  If $G,H$ are cographs, their disjoint union $G\cup H$ is a cograph.

\end{enumerate} 
\end{definition}

By other words, cographs are the smallest class of finite graph, that contain all one vertex graph, and it is closed for operations 2,3. There are many cographs, we can construct. Empty graphs are cographs, as they are the disjoint union of their vertices as one vertex graphs. Complete graphs are the complement of empty graphs, so they are also cographs. Complete multipartite graphs are also cographs, as they are the complement of disjoint union of complete graphs, and so on. 
Cographs can be characterized in different ways, given by  \cite[Theorem 2]{CLS}.   


\begin{theorem}\label{cograph_characterization}
 For a finite graph $G$ the following statements are equivalent:
 \begin{enumerate}[1)]
 \item   $G$ is a cograph.
\item  For any $U\subseteq V[G]$ and $X,Y\subseteq U$, if $X$ is a maximal clique in $G[U]$ and $Y$ is a maximal independent set of $G[U]$, then $X\cap Y\neq \emptyset $.
 
\item  $G$ is $P_{4}$-free (it does not contain a path of 4 vertices as an induced subgraph).

 \end{enumerate}
\end{theorem}

Now we characterize the trivially power-colorable cographs.

\begin{lemma}\label{cliqued_cographs}
 Let $G$ be a cograph, with $\chi (G)=k$, then the following statements are equvivalent:
\begin{enumerate}[1)]
\item  $G$ is tight.
 \item  There is a coloring $\phi :V(G)\rightarrow k$, that is tight.
 \item  For all maximal cliques of $X\subseteq V(G)$, we have $|X|=k$.
 \item $G$ is strongly-cliqued.
 \item $G$ is weakly-cliqued.
\end{enumerate} 
 
\end{lemma}

\begin{proof}
$(1)\rightarrow (2)$: Obvious

$(2)\rightarrow (3)$: For all $0\le l<k$, $Y_{l}=\{ v\in V(G):\phi (v)=l\} $ is an independent set. If $u\not\in Y_{l}$, then $\phi (u)\neq l$, and since $\phi $ is tight, there is some $w\in V(G)$ with $uw\in E(G)$ and $\phi (w)=l$. Then $w\in Y_{l}$, so $Y_{l}\cup \{ u\} $ is not independent, thus $Y_{l}$ is maximal independent for all $l$. If $X\subseteq V(G)$ is a maximal clique, then by (2) of Theorem \ref{cograph_characterization}, we have $X\cup Y_{l}\neq \emptyset$ for all $l$. Since $X$ intersects the pairwise disjoint sets $Y_{0},...,Y_{k-1}$, we have $|X|\ge k$. On the other hand $\chi (G)=k$, so $G$ can not have any clique larger than $k$, thus $|X|=k$.

$(3)\rightarrow (4)$: If $v\in V(G)$ were an isolated vertex, then $\{ v\} $ would be itself a maximal clique of size less then $k$, so $G$ does not have any isolated vertices. If $u,v\in V(G)$ with $uv\in E(G)$, then $\{ u,v\} \subseteq V(G)$ is a clique, so there is a maximal clique $X\subseteq V(G)$, with $\{ u,v\} \subseteq X$. Then $|X|=k$ and $u,v\in X$, so $G$ is strongly-cliqued.

$(4)\rightarrow (5)\rightarrow (1)$: We already had for general graphs.
\end{proof}

\begin{theorem}\label{trivially_colorable_cographs}
 A cograph $G$ is trivially power-colorable if and only if it is tight, connected and $\chi (G)\ge 3$.
\end{theorem}

\begin{proof}
 By Theorem \ref{trivial_connected}, Theorem \ref{trivial_chromatic_3} and Theorem \ref{trivial_tight}, these conditions are nessesary. If a cograph $G$ has these properties, then since $G$ is tight, by Lemma \ref{cliqued_cographs}, it is weakly-cliqued, and then by Theorem \ref{weakly_cliqued_power_trivial} it is trivially power-colorable.
\end{proof}

\section{Colorings of infinite products}

Now we turn our attention to the colorings of $G^{\lambda }$, where $G$ is a 
trivially power-colorable finite graph and $\lambda $ is an infinite cardinal. 

As we will see in Observation \ref{obs:color},  
 we cannot expect a single coordinate to determine 
the color of a vertex in an arbitrary coloring of $G^{\lambda}$.
However, in  constrast to this result,  Theorem \ref{index_ultrafilter} provides  a 
complete structural 
description of  any proper coloring of $G^{\lambda}$. 

Fist, we define a natural class of colorings of $G^{\lambda}$.

\begin{definition}\label{df:uf-color}
Let $G$ be a finite graph,  $\phi:V(G)\to k$  a coloring, and 
$\mathcal  U$ be  an ultrafilter on an infinite cardinal ${\lambda}$.
Define the coloring $\phi^{\mathcal U}:V(G^{\lambda})\to k$ as follows:
for each  $\mathbf v=(v_{\alpha})_{{\alpha}<{\lambda}}\in V(G^{\lambda})$
pick the unique $\mathbf v_{\mathcal U}\in V(G)$ such that
\begin{displaymath}
\{{\alpha}<{\lambda}: v_{\alpha}=\mathbf v_{\mathcal U}\}\in \mathcal U,
\end{displaymath} 
and let $\phi^{\mathcal U}(\mathbf v)=\phi(\mathbf v_{\mathcal U})$. 
(Since $V(G)$ is finite, $\mathbf v_{\mathcal U}$ is defined.)
\end{definition}

\begin{observation}\label{obs:color}
(1) If $G$ is  a finite graph,  $\phi:V(G)\to k$ is a proper coloring of $G$,  and 
 $\mathcal  U$ is an ultrafilter on an infinite cardinal ${\lambda}$, then 
$\phi^{\mathcal U}:V(G^{\lambda})\to k$ is a proper coloring of $G^{\lambda}$.

\noindent (2)  If $\mathcal U$ is non-principal, then finitely many coordinates
do not  determine  the $\phi^{\mathcal U}$-color  of a vertex of $G^{\lambda}$. In particular,
$\phi^{\mathcal U}$ is not a trivial coloring. 
\end{observation}

The first part of following observations are straightforward, based on that fact that  
if $\mathbf v\mathbf w$ is an edge in $G^{\lambda}$, then 
$\mathbf v_{\mathcal U}\mathbf w_{\mathcal U}$ is an edge in $G$.
The second part follows from the fact that if 
$\{{\alpha}<{\lambda}:\mathbf v({\alpha})=\mathbf w({\alpha})\}\in \mathcal U$,
then $\phi^{\mathcal U}(\mathbf v)=\phi^{\mathcal U}(\mathbf w)$.

Komjáth and Totik in \cite{KoTo} has shown, that for all $k\ge 3$, the colorings of the graph $K_{k}^{\lambda }$ are in this form, defined by an ultrafilter. Now we generalize this theorem to all trivially power-colorable finite graphs.

\begin{theorem}\label{trivial_infinite_product}
  If $G$ is a trivially power-colorable finite graph, and  $\Phi $ is a proper $\chi(G)$-coloring of $G^{\lambda }$ for some infinite cardinal $\lambda $, then there is an ultrafilter $\mathscr{U}$  on $\lambda $ and a $\chi(G)$-coloring $\phi $ on $G$ such that $\Phi=\phi^{\mathscr{U}}$.
 \end{theorem}
 
 \begin{proof}
Write $k=\chi(G)$. 
  The coloring   $\Phi:V(G^{\lambda})\to k$ is fixed.  


\begin{definition}\label{partition}
 A {\em finite partition of ${\lambda}$} is a  partition 
 of ${\lambda}$ into finitely many non-empty pieces. 
 Let $Part(\lambda )$ be the set of all finite partitions of $\lambda $. 
 For a finite partition $\mathcal{P}$, let 
 $$V_{\mathcal{P}}=\{ v\in V(G^{\lambda }):v\restriction_{A}\text{ is constant for all } A\in\mathcal{P}\}$$ and 
 $G_{\mathcal{P}}=G^{\lambda}[V_{\mathcal{P}}]$.
\end{definition}

If $\mathbf v=(v_{\xi})_{{\xi}<{\lambda}}\in V_{\mathcal{P}}$, and $A\in \mathcal{P}$, then, by definition, there is a $\mathbf v_{A}\in V(G)$, such that $v_{\xi }=v_{A}$ for all $\xi \in A$. 
 If $\mathcal P=\{A_1,\dots, A_n\}$, then the map
 \begin{displaymath}
 \mathbf v \mapsto \langle \mathbf v_{A_1},\dots, \mathbf v_{A_n}\rangle
 \end{displaymath}
 gives an isomorphism between $G_{\mathcal{P}}$ and $G^n$.

Now  
 $\Phi |_{V_{\mathcal{P}}}$ is a $k$-coloring of $G_{\mathcal{P}}$. Since $G$ is trivially power-colorable, and $G_{\mathcal{P}}$ and $G^n$ are isomorphic,   there is some 
 $1\le i\le n$, 
 and a coloring $\phi$ of $G$, such that 
 $$\Phi (\mathbf v)=\phi(\mathbf v_{A_{i}})$$ for all $\mathbf v\in V_{\mathcal{P}}$. 
 Write $i_{\mathcal P}=i$, $\phi_{\mathcal P}=\phi$ and  
 $I(\mathcal{P})=A_{i}$. 
 

\begin{definition}\label{refinement}
 For $\mathcal{P}_{1},\mathcal{P}_{2}\in Part(\lambda )$, we say, that 
 $\mathcal{P}_{1}$ is a {\em refinement} of $\mathcal{P}_{2}$, denoted by 
 $\mathcal{P}_{1}\prec \mathcal{P}_{2}$ if for all $A\in \mathcal{P}_{1}$, there is a $B\in \mathcal{P}_{2}$ with $A\subseteq B$. 
 
 For $\mathcal{P}_{1},\mathcal{P}_{2}\in Part(\lambda )$, their 
 {\em greatest common refinement} is $$\mathcal{P}_{1}\sqcap \mathcal{P}_{2}=\{ A\cap B:A\in \mathcal{P}_{1}, B\in \mathcal{P}_{2},A\cap B\neq \emptyset \}.$$ It is easy to see, that this is also a finite partition, that is the refinement of both $\mathcal{P}_{1}$ and $\mathcal{P}_{2}$. 
\end{definition}

If $\mathcal{P}_{1}\prec \mathcal{P}_{2}$, then $V_{\mathcal{P}_{2}}\subseteq V_{P_{1}}$. For the trivial partition $\mathcal{P}=\{\lambda \}$, $V_{\mathcal{P}}$ contains only the constant functions.

\begin{lemma}\label{refinement_index}
 If $\mathcal{P}_{1},\mathcal{P}_{2}\in Part(\lambda )$ with $\mathcal{P}_{1}\prec \mathcal{P}_{2}$, then $I(\mathcal{P}_{1})\subseteq I(\mathcal{P}_{2})$.
\end{lemma}

\begin{proof}
Assume on the contrary that $I(\mathcal{P}_{1})\cap I(\mathcal{P}_{2})=\emptyset.$
Pick $B\in \mathcal P_2$ with $I(\mathcal{P}_{1})\subset B$. Let $vw$ be an edge in 
$G$. Consider the vertices $\mathbf u, \mathbf v$ and  $\mathbf w$ from $V(G^{\lambda})$,
where $\mathbf v({\xi})=v$ and $\mathbf w({\xi})=v$ for each ${\zeta}<{\lambda}$, and 
\begin{displaymath}
{\mathbf{u}({\zeta})}=\left\{\begin{array}{ll}
{v}&\text{if ${\xi}\in B$ }\\
{w}&\text{otherwise. }\\
\end{array}\right.
\end{displaymath} 
Since $\mathbf u\in V_{\mathcal P_2}$ and $\mathbf w\restriction I(\mathcal P_2)=\mathbf u\restriction I(\mathcal P_2)$, 
we have $\Phi(\mathbf w)=\Phi(\mathbf u)$.
Since $\mathbf u\in V_{\mathcal P_1}$ and $\mathbf v\restriction I(\mathcal P_1)=\mathbf u\restriction I(\mathcal P_1)$, 
we have $\Phi(\mathbf w)=\Phi(\mathbf u)$. But $\mathbf v\mathbf w$ is en edge in $G^{\lambda}$. Contradiction.  
\end{proof}


\begin{lemma}\label{common_refinement_index}
 If $\mathcal{P}_{1},\mathcal{P}_{2}\in Part(\lambda )$, then $I(\mathcal{P}_{1}\sqcap \mathcal{P}_{2})=I(\mathcal{P}_{1})\cap I(\mathcal{P}_{2})$.
\end{lemma}

\begin{proof}
 Since $\mathcal{P}_{1}\sqcap \mathcal{P}_{2}\prec \mathcal{P}_{1}$ and $\mathcal{P}_{1}\sqcap \mathcal{P}_{2}\prec \mathcal{P}_{2}$, by Lemma \ref{refinement_index}, we have $I(\mathcal{P}_{1}\sqcap \mathcal{P}_{2})\subseteq I(\mathcal{P}_{1})$ and $I(\mathcal{P}\sqcap \mathcal{P}_{2})\subseteq I(\mathcal{P}_{2})$, so $I(\mathcal{P}_{1}\sqcap \mathcal{P}_{2})\subseteq I(\mathcal{P}_{1})\cap I(\mathcal{P}_{2})$. This means, that $I(\mathcal{P}_{1})\cap I(\mathcal{P}_{2})\neq \emptyset $ and thus $I(\mathcal{P}_{1})\cap I(\mathcal{P}_{2})\in \mathcal{P}_{1}\sqcap \mathcal{P}_{2}$. Then $I(\mathcal{P}_{1}\sqcap \mathcal{P}_{2})$ is an element of $\mathcal{P}_{1}\sqcap \mathcal{P}_{2}$, contained in $I(\mathcal{P}_{1})\cap I(\mathcal{P}_{2})$, but $\mathcal{P}_{1}\sqcap \mathcal{P}_{2}$ is a finite partition, so $I(\mathcal{P}_{1}\sqcap P_{2})=I(\mathcal{P}_{1})\cap I(\mathcal{P}_{2})$.
\end{proof}

\begin{lemma}\label{index_ultrafilter}
 The set $\mathscr{U}=\{ I(\mathcal{P}):\mathcal{P}\in Part (\lambda )\}$ is an ultrafilter on $\lambda $.
\end{lemma}

\begin{proof}
 Clearly $\emptyset \not\in \mathscr{U}$, as it is not even an element of a partition. For the trivial partition, we have $\lambda =I(\{ \lambda\} )\in \mathscr{U}$. 
 If $A,B\in \mathscr{U}$, then there are some $\mathcal{P}_{1},\mathcal{P}_{2}\in Part(\lambda )$ with $A=I(\mathcal{P}_{1}),B=I(\mathcal{P}_{2})$. 
 Then by Lemma \ref{common_refinement_index}, we have $A\cap B=I(\mathcal{P}_{1})\cap I(\mathcal{P}_{2})=I(\mathcal{P}_{1}\sqcap \mathcal{P}_{2})\in \mathscr{U}$. 
 Now let $A\in \mathscr{U}$ and $B\subseteq \lambda $ with $A\subseteq B$. Let $\mathcal{P}_{1}\in Part(\lambda )$ be such 
 that $A=I(\mathcal{P}_{1})$ and let $\mathcal{P}_{2}\in Part(\lambda )$ be any finite partition with $B\in \mathcal{P}_{2}$ (for example $\mathcal{P}_{2}=\{ B,\lambda \setminus B\}$). 
 Then $A=A\cap B\in \mathcal{P}_{1}\sqcap \mathcal{P}_{2}$, and by Lemma \ref{refinement_index}, since $\mathcal{P}_{1}\sqcap \mathcal{P}_{2}\prec \mathcal{P}_{1}$, we have $I(\mathcal{P}_{1}\sqcap \mathcal{P}_{2})\subseteq I(\mathcal{P}_{1})=A$, and since $\mathcal{P}_{1}\sqcap \mathcal{P}_{2}$ is a partition, we must have $I(\mathcal{P}_{1}\sqcap \mathcal{P}_{2})=A$. 
 On the other hand, we have $\mathcal{P}_{1}\sqcap \mathcal{P}_{2}\prec \mathcal{P}_{2}$, by Lemma \ref{refinement_index}, 
 we have $A=I(\mathcal{P}_{1}\sqcap \mathcal{P}_{2})\subseteq I(\mathcal{P}_{2})$. 
 But since $\mathcal{P}_{2}$ is also a partition, $B$ is its only element containing $A$, so $B=I(\mathcal{P}_{2})\in \mathscr{U}$. Finally, if $A\subseteq \lambda $, we need to show that either $A$ or $\lambda \setminus A$ is in $\mathscr{U}$. It is obvious, if one of them is $\emptyset$, so suppose $A,\lambda \setminus A\neq \emptyset $. 
 Then $\mathcal{P}=\{A,\lambda \setminus A\}$ is a finite partition, so $I(\mathcal{P} )\in \mathscr{U}$ is either $A$ or $\lambda \setminus A$, and that is what we needed to prove.
\end{proof}

We are ready to conclude  the Proof of Theorem \ref{trivial_infinite_product}.
 Let $\mathscr{U}$ be the ultrafilter, defined in Lemma \ref{index_ultrafilter}. 
 First we will show, that if $u,v\in V(G^{\lambda })$ with $\{ \xi \in \lambda :u_{\xi}=v_{\xi}\} \in \mathscr{U}$, then $\Phi (u)=\Phi (v)$. 
 Let $$\mathcal{P}(u)=\{ \{ \xi \in \lambda :u_{\xi }=w\} :w\in V(G)\} \setminus\{ \emptyset\}\in Part(\lambda )$$ and similarly,  $$\mathcal{P}(v)=\{ \{ \xi \in \lambda :v_{\xi }=w\} :w\in V(G)\} \setminus\{ \emptyset\}\in Part(\lambda ).$$ Then $u\in V_{\mathcal{P}(u)}, v\in V_{\mathcal{P}(v)}$, so $u,v\in V_{\mathcal{P}(u)\sqcap \mathcal{P}(v)}$, so they are both constant on $I(\mathcal{P}(u)\sqcap \mathcal{P}(v))$. This way it is either $u_{\xi }=v_{\xi }$ for all $\xi \in I(\mathcal{P}(u)\sqcap \mathcal{P}(v))$, or $u_{\xi }\neq v_{\xi }$ for all $\xi \in I(\mathcal{P}(u)\sqcap \mathcal{P}(v))$. On the other hand, we have $I(\mathcal{P}(u)\sqcap \mathcal{P}(v))\in \mathscr{U}$, so $I(\mathcal{P}(u)\sqcap \mathcal{P}(v))\cap \{ \xi \in \lambda :u_{\xi}=v_{\xi}\} \neq \emptyset $, so the second one is not an option. Thus we have $u_{I(\mathcal{P}(u)\sqcap \mathcal{P}(v))}=v_{I(\mathcal{P}(u)\sqcap \mathcal{P}(v))}$, so $\Phi (u)=\Phi (v)$.
 
 Now for any $w\in V(G)$, let $c(w)\in V(G^{\lambda })$, be such that $c(w)_{\xi }=w$ for all $\xi \in \lambda $. Let as define a coloring $\phi $ on $G$, such that $\phi (w)=\Phi (c(w))$ for $w\in V(G)$. We need to show, that this is a good coloring of $G$. If $w,w'\in V(G)$, with $ww'\in E(G)$, then $c(w)c(w')\in E(G^{\lambda })$, so $\phi (w)=\Phi (c(w))\neq \Phi (c(w'))=\phi (w')$. Finally let $v\in V(G^{\lambda })$ be arbitrary. Then $\{ \xi \in \lambda :v_{\xi}=c(v_{\mathscr{U}})_{\xi}\} =\{ \xi \in \lambda :v_{\xi}=v_{\mathscr{U}}\} \in \mathscr{U}$, so we have $\Phi (v)=\Phi(c(v_{\mathscr{U}}))=\phi (v_{\mathscr{U}})$, and then we are done.
\end{proof}

\section{Conclusion and problems}

We provided a partial characterization of trivially power-colorable finite graphs.
 By Theorems \ref{trivial_connected}, \ref{trivial_chromatic_3} and \ref{trivial_tight}, a trivially power-colorable graph $G$ must be connected, tight and satisfy $\chi (G)\ge 3$. COnversely,  by Theorem \ref{weakly_cliqued_power_trivial}, if $G$ is a weakly cliqued connected graph with $\chi (G)\ge 3$, then $G$ is trivially power-colorable.

These conditions do not always align, as there are examples of finite graphs that are tight, but not weakly cliqued. However, we do not yet know whether these graphs are trivially power colorable, or not, which raises two open problems,  one of which must have positive solution.

\begin{problem}
    Is there any finite graph $G$ that is trivially power-colorable but not weakly cliqued?
\end{problem}

\begin{problem}
    Is there any finite graph $G$, that is connected, tight, $\chi (G)\ge 3$ holds, but is not trivially power-colorable?
\end{problem}

For cographs, the nessecary an sufficient conditions align,  characterizing  trivially power-colorability. This raises the question of whether  this result can be extended to a larger class of graphs.

\begin{problem}
    Is there a "natural" class $\Gamma $ of finite graphs larger than the class of cographs, such that one of the following is always holds?

    - If $G\in \Gamma $ is connected, tight and $\chi (G)\ge 3$, then $G$ is weakly cliqued.

    - If $G\in \Gamma $ is connected, tight and $\chi (G)\ge 3$, then $G$ is trivially power-colorable.

    - If $G\in \Gamma$ is trivially power-colorable, then $G$ is weakly cliqued.
\end{problem}

By Theorem \ref{trivial_infinite_product}, we also know that the colorings of an infinite power of a finite, trivially power-colorable graphs are determined by an ultrafilter. However,   we have not yet addressed  the case of infinite graphs. If a graph $G$ has infinite chromatic number, then it is easy to see that they cannot be trivially power-colorable, but in case $\chi (G)$ is finite, trivial power colorability can be defined in a similar way. This leads to the  following problem.

\begin{problem}
    Is it true that if $G$ is an infinite graph with finite chromatic number $\chi(G)=k$, such that for all $n$, every $k$-colorings of $G^{n}$ is trivial, then for any infinite $\lambda $, all colorings of $G^{\lambda }$ are determined by an ultrafilter?
\end{problem}

\bibliographystyle{plain}

\end{document}